\documentclass[12pt]{amsart}
\usepackage{amscd,amssymb,amsmath}
\newtheorem{theorem}{Theorem}[section]

\theoremstyle{definition}

\newtheorem{proposition}[theorem]{Proposition}
\newtheorem{corollary}[theorem]{Corollary}

\date \kill
\theoremstyle{remark}
\newtheorem{remark}{\bf Remark}[section]

\numberwithin{equation}{section}
\begin{document}
\title{A note on generalized equivariant homotopy groups}
\author{Marek Golasi\'nski}
\address{Faculty of Mathematics and Computer Science, Nicolaus Copernicus University, Chopina 12/18, 87-100, Toru\'n, Poland\newline
Faculty of Mathematics and Computer Science, University of Warmia and Mazury,
\.Zo\l nierska 14, 10-561 Olsztyn, Poland}
\email{marek@mat.uni.torun.pl\newline marek@matman.uwm.edu.pl}
\author{Daciberg Gon\c calves}
\address{Dept. de Matem\'atica - IME - USP, Caixa Postal 66.281 - CEP 05311-970,
S\~ao Paulo - SP, Brasil}
\email{dlgoncal@ime.usp.br}
\author{Peter Wong}
\address{Department of Mathematics, Bates College, Lewiston,
ME 04240, U.S.A.}
\email{pwong@bates.edu}
\thanks{
This work was conducted during the second and third authors' visit to the Faculty
of Mathematics and Computer Science, Nicolaus Copernicus University July 27 - 31, 2007.
The second and third authors would like to thank the Faculty of Mathematics and
Computer Science for its hospitality and support}

\begin{abstract}
In this paper, we generalize the equivariant homotopy groups or equivalently the
Rhodes groups. We establish a short exact sequence relating the generalized Rhodes
groups and the generalized Fox homotopy groups and we introduce $\Gamma$-Rhodes groups,
where $\Gamma$ admits a certain co-grouplike structure.  Evaluation subgroups of $\Gamma$-Rhodes groups are discussed.
\end{abstract}
\date{\today}
\maketitle

\newcommand{\af}{\alpha}
\newcommand{\et}{\eta}
\newcommand{\ga}{\gamma}
\newcommand{\ta}{\tau}
\newcommand{\ph}{\varphi}
\newcommand{\bt}{\beta}
\newcommand{\lb}{\lambda}
\newcommand{\wh}{\widehat}
\newcommand{\wt}{\widetilde}
\newcommand{\sg}{\sigma}
\newcommand{\om}{\omega}
\newcommand{\cH}{\mathcal H}
\newcommand{\cF}{\mathcal F}
\newcommand{\N}{\mathcal N}
\newcommand{\R}{\mathcal R}
\newcommand{\Ga}{\Gamma}
\newcommand{\cc}{\mathcal C}

\newcommand{\bea} {\begin{eqnarray*}}
\newcommand{\beq} {\begin{equation}}
\newcommand{\bey} {\begin{eqnarray}}
\newcommand{\eea} {\end{eqnarray*}}
\newcommand{\eeq} {\end{equation}}
\newcommand{\eey} {\end{eqnarray}}

\newcommand{\ovl}{\overline}
\newcommand{\vv}{\vspace{4mm}}
\newcommand{\lra}{\longrightarrow}

\bibliographystyle{plain}



\section*{Introduction}

In 1966, F. Rhodes \cite{rhodes1} introduced the fundamental group of a transformation group $(X,G)$ for a topological space on which a group $G$ acts. This group, denoted by $\sigma_1(X,x_0,G)$, is the equivariant analog of the classical fundamental group $\pi_1(X,x_0)$. Rhodes showed that $\sigma_1(X,x_0,G)$ is a group extension of $\pi_1(X,x_0)$ with quotient $G$. Thus, $\sigma_1(X,x_0,G)$ incorporates the $G$-action as well as the action of $\pi_1(X,x_0)$ on the universal
cover $\tilde X$ of the space $X$. This group has been used in \cite{wong} to study the Nielsen fixed point theory for equivariant maps. In 1969, F. Rhodes \cite{rhodes2} extended $\sigma_1(X,x_0,G)$ to $\sigma_n(X,x_0,G)$, which is the equivariant higher homotopy group of $(X,G)$. Like $\sigma_1(X,x_0,G)$, $\sigma_n(X,x_0,G)$ is an extension of the Fox torus homotopy group $\tau_n(X,x_0)$ but not of the classical homotopy group $\pi_n(X,x_0)$ by $G$. The Fox torus homotopy groups were first introduced by R. Fox \cite{fox} in 1948 in order to give a geometric interpretation of the classical Whitehead product. Recently, a modern treatment of $\tau_n(X,x_0)$ and of $\sigma_n(X,x_0,G)$ has been given in \cite{ggw} and in \cite{ggw2}, respectively.  In \cite{ggw2}, we further investigated the relationships between the Gottlieb groups of a space and of its orbit space, analogous to the similar study in \cite{gg2}. Further properties of the Fox torus homotopy groups, their generalizations, and Jacobi identities were studied in \cite{ggw3}. It is therefore natural to generalize $\sigma_n(X,x_0,G)$ to more general constructions with respect to general spaces and to co-grouplike spaces $\Gamma$
other than the $1$-sphere $\mathbb{S}^1$.

The main objective of this paper is to generalize $\sigma_n(X,x_0,G)$ of a $G$-space $X$ with respect to a space $W$ and also
with respect to a pair $(W,\Gamma)$, where $W$ is a space and $\Gamma$ satisfies a suitable notion of the
classical co-grouplike space. We prove in section 1 that the Rhodes exact sequence of \cite{rhodes2} can be
generalized to $\sigma_W(X,x_0,G):=\{[f;g]|f:(\wh \Sigma W, v_1,v_2) \to (X, x_0, gx_0)\}$,
the $W$-Rhodes group, with the generalized Fox torus homotopy group
$\tau_W(X,x_0)$ as the kernel. In section 2, we further extend the construction of Rhodes groups to
$\sigma^{\Gamma}_W(X,x_0,G):=\{[f;g]|\;f:(\Gamma(W),\bar{\gamma}_1,\bar{\gamma}_2) \to (X,x_0,gx_0)\}$,
the $W$-$\Gamma$-Rhodes groups, where $\Gamma$ admits a co-grouplike structure with {\it two} basepoints.
Under such assumptions, 
we obtain a $W$-$\Gamma$-generalization of the Rhodes exact sequence \cite{rhodes2}.
In the last section, we generalize the notion of the Gottlieb (evaluation) subgroup to that of a $W$-$\Gamma$-Rhodes group and we establish a short exact sequence generalizing \cite[Theorem 2.2]{ggw2}.
Throughout, $G$ denotes a group acting on a compactly generated Hausdorff path-connected space $X$
with a basepoint $x_0$.
The associated pair $(X,G)$ is called in the literature a transformation group.

\section{Generalized Rhodes groups}

For $n\ge 1$, F.\ Rhodes \cite{rhodes2} defined higher homotopy groups $\sigma_n(X,x_0,G)$ of a pair $(X,G)$
which is an extension of $\tau_n(X,x_0)$ by $G$ so that
\begin{equation}\label{rhodes-exact}
1\to \tau_n(X,x_0) \to\sigma_n (X,x_0,G) \to G \to 1
\end{equation}
is exact. Here, $\tau_n(X,x_0)$ denotes the $n$-th torus homotopy group of $X$ introduced
by R.\ Fox \cite{fox}. The group $\tau_n=\tau_n(X,x_0)$ is defined to be the fundamental group
of the function space $X^{\mathbb{T}^{n-1}}$ and is uniquely determined by the groups $\tau_1, \tau_2, \ldots, \tau_{n-1}$
and the Whitehead products, where $\mathbb{T}^{n-1}$ is the $(n-1)$-dimensional torus. The group $\tau_n$ is
non-abelian in general.

Now we recall the construction of $\sigma_n(X,x_0,G)$ presented in \cite{rhodes2}. Suppose that $X$ is a $G$-space with a basepoint $x_0\in X$ and
let $C_n=I\times \mathbb{T}^{n-1}$. We say that a map $f : C_n\to X$ is of {\em order}
$g\in G$ provided $f(0,t_2,\ldots,t_n)=x_0$ and $f(1,t_2,\ldots,t_n)=g(x_0)$ for $(t_2,\ldots,t_n)\in \mathbb{T}^{n-1}$.
Two maps $f_0,f_1: C^n\to X$ of order $g$ are said to be homotopic if there exists a continuous map $F : C^n\times I\to X$ such that:
\begin{itemize}
\item $F(t,t_2,\ldots,t_n,0)=f_0(t,t_2,\ldots,t_n)$;
\item $F(t,t_2,\ldots,t_n,1)=f_1(t,t_2,\ldots,t_n)$;
\item $F(0,t_2,\ldots,t_n,s)=x_0$;
\item $F(1,t_2,\ldots,t_n,s)=gx_0$
for all $(t_2,\ldots,t_n)\in \mathbb{T}^{n-1}$ and $s,t\in $I.
\end{itemize}
Denote by $[f;g]$ the homotopy class of a map $f: C_n\to X$ of order $g$ and by $\sigma_n(X,x_0,G)$ the set of all such homotopy classes.
We define an operation $*$ on the set $\sigma_n(X,x_0,G)$ by
$$[f';g']*[f;g]:=[f' + g'f;g'g].$$
This operation makes $\sigma_n(X,x_0,G)$ a group.

We have generalized in \cite{ggw} the Fox torus homotopy groups. In this section, we give a similar
generalization of Rhodes groups. In a special case, we obtain an extension group of the Abe group
considered in \cite{abe}.

Let $X$ be a path-connected space with a basepoint $x_0$. For any space $W$, we let
$$
\sigma_W(X,x_0,G):=\{[f;g]|f:(\wh \Sigma W, v_1,v_2) \to (X, x_0, gx_0)\}
$$
where $[f;g]$ denotes the homotopy class of the map $f$ of order $g\in G$, $v_1$ and $v_2$ are
the vertices of the cones $C^+W$ and $C^-W$, respectively and $\wh \Sigma W=C^+W \cup C^-W$.
Under the operation $[f_1;g_1]*[f_2;g_2]:=[f_1 + g_1f_2;g_1g_2]$, $\sigma_W$ is a group called a $W$-{\em Rhodes group}.
\par Write $C(W,X)$ for the mapping space of all continuous maps from $W$ to $X$ with the compact-open topology.
We point out that $\sigma_W(X,x_0,G)=\sigma_1(C(W,X),\bar{x_0},G)$ provided $W$ is a locally-compact space, where $(gf)(x)=gf(x)$ for
$f\in C(W,X)$, $g\in G$ and $\bar{x_0}$ denotes the constant map from $C(W,X)$ determined by the point $x_0\in X$.
\par The canonical projection $\sigma_W(X,x_0,G)\to G$ given by $[f;g]\mapsto g$ for $[f;g]\in\sigma_W(X,x_0,G)$
has the kernel $\{[f;1]|f:(\wh \Sigma W, v_1,v_2) \to (X, x_0, x_0)\}$. It is easy to see that this kernel
is isomorphic to the generalized Fox torus group $[\Sigma (W\sqcup *),X]=\tau_W(X,x_0)$ defined in \cite{ggw}.
Therefore, we get the following result.

\begin{theorem}\label{general-exact}
The following sequence
\begin{equation}\label{generalized-rhodes-exact}
1 \to \tau_W(X,x_0) \to \sigma_W(X,x_0,G) \to G \to 1
\end{equation}
is exact.
\end{theorem}

\begin{remark} When $W=\mathbb{T}^{n-1}$, the $(n-1)$-dimensional torus, $\sigma_W$ coincides with the $n$-th Rhodes group $\sigma_n$ and \eqref{generalized-rhodes-exact} reduces to \eqref{rhodes-exact}.
When $W=\mathbb{S}^{n-1}$, the $(n-1)$-sphere, $\tau_W$ becomes $\kappa_n$, the $n$-th Abe group (see \cite{fox} or \cite{ggw}). Thus, by Theorem \ref{general-exact}, we have the following exact sequence
\begin{equation}\label{abe}
1 \to \pi_n(X,x_0)\rtimes \pi_1(X,x_0)\cong \kappa_n(X,x_0) \to \sigma_{\mathbb{S}^{n-1}}(X,x_0,G) \to G \to 1.
\end{equation}
\end{remark}

One can also generalize the split exact sequence for Rhodes groups from \cite{rhodes2} as follows.

\begin{theorem}\label{general-rhodes-fox-split}
Let $W$ be a space with a basepoint $w_0$. Then, for any space $V$, the following sequence
\begin{equation}\label{general-rhodes-split}
1\to [(V\times W)/V,\Omega X] \to \sigma_{V\times W}(X,x_0,G) \stackrel{\dashleftarrow}{\to} \sigma_V(X,x_0,G) \to 1
\end{equation}
is split exact.
\end{theorem}
\begin{proof} By \cite[Theorem 3.1]{ggw}, we have the following split
exact sequence
\begin{equation}\label{ggw-sequence}
1\to [(V\times W)/V,\Omega X] \to \tau_{V\times W}(X) \stackrel{\dashleftarrow}{\to} \tau_V(X) \to 1.
\end{equation}
Given $[F;g]\in \sigma_{V\times W}(X,x_0,G)$, where $F:\wh \Sigma (V\times W) \to X$,
let $f:\wh \Sigma V \to X$ be the composite map of $\wh \Sigma V \approx \wh \Sigma (V\times \{w_0\}) \to \wh \Sigma (V\times W)$
with $F$. This map gives rise to a homomorphism $\sigma_{V\times W}(X,x_0,G) \to \sigma_V(X,x_0,G)$.
Likewise, using the projection $V\times W\to V$, one obtains a section
$\sigma_V(X,x_0,G) \to \sigma_{V\times W}(X,x_0,G)$. We have the
following commutative diagram
\begin{equation}\label{general-rhodes-fox-split-sequence}
\begin{CD}
1 @>>>    \tau_{V\times W}(X,x_0)       @>>>  \sigma_{V\times W}(X,x_0,G)     @>>>   G    @>>> 1\\
@.        @VVV                 @VVV                   @|         \\
1 @>>>    \tau_{V}(X,x_0)   @>>>  \sigma_{V}(X,x_0,G) @>>>   G    @>>> 1,
\end{CD}
\end{equation}
where the first two vertical homomorphisms have sections.
Combining with \eqref{ggw-sequence}, the assertion follows.
\end{proof}

As an immediate corollary of Theorem \ref{general-rhodes-fox-split}, we have the following:

\begin{corollary}\label{general-rhodes-corollary}
{\em The following sequence
\begin{equation}\label{general-rhodes-corollary-split}
1\to [W,\Omega X] \to \sigma_W(X,x_0,G) \stackrel{\dashleftarrow}{\to} \sigma_1(X,x_0,G) \to 1
\end{equation}
is split exact.}
\end{corollary}
\begin{proof} The result follows from Theorem \ref{general-rhodes-fox-split} by letting $V$ be a point.
\end{proof}

\begin{remark} For any space $W$, Corollary \ref{general-rhodes-corollary} asserts that $\sigma_1(X,x_0,G)$
acts on $[\Sigma W,X]=[W,\Omega X]$ according to the splitting. Furthermore, when $W=\mathbb{S}^{n-1}$, this corollary gives an
alternate description of the action of $\sigma_1$ on $\pi_n(X)$ as described in \cite[Remark 1.4]{ggw2}.
In this case, $\sigma_W(X,x_0,G)=\sigma_{\mathbb{S}^{n-1}}(X,x_0,G)$ is the extension group of the $n$-th Abe
group $\kappa_n(X,x_0)$ \cite{abe} as in \eqref{abe}. Thus, one can either embed $\sigma_1$ in $\sigma_n$ as in
\cite[Remark 1.4]{ggw2} or in $\sigma_{\mathbb{S}^{n-1}}(X,x_0,G)$.
\end{remark}

Unlike the reduced suspension $\Sigma$ which has the loop functor
$\Omega$ as its right adjoint, the un-reduced suspension $\wh \Sigma$
does {\it not} admit a right adjoint. Nevertheless, one can describe the
adjoint property for the $W$-Rhodes groups as follows.
Recall that a typical element in $\sigma_W(X,x_0,G)$ is a homotopy class
$[f;g]$ where $f:(\wh\Sigma W, v_1,v_2) \to (X,x_0,gx_0)$. Thus,
$\sigma_W$ is a subset of $[\wh \Sigma W,X]_0 \times G$, where
$[\wh \Sigma W,X]_0$ denotes the homotopy classes of maps
$f:\wh \Sigma W \to X$ such that $f(v_1)=x_0$ and $f(v_2)$ is independent
of the homotopy class of $f$.
Then, $\sigma_W$ is also a subset of $[W,\mathcal P_{x_0}]^* \times G$,
where $[W,\mathcal P_{x_0}]^*$ denotes the set of homotopy classes of unpointed maps from $W$ to
the space $\mathcal P_{x_0}$ of paths originating from $x_0$. In the special case when
$G=\{1\}$, $\sigma_W=[\Sigma (W\cup *),X]=\sigma_W^*=[W,\Omega Y]^*=[W\cup *,\Omega X]$.

\section{Generalized $W$-$\Gamma$-Rhodes groups}

In the definition of the generalized Rhodes group $\sigma_W(X,x_0,G)$, the two
cone points from the un-reduced suspension $\wh \Sigma W=C^+W\cup C^-W$
play an important role. Therefore in replacing $\mathbb{S}^1$ with arbitrary co-grouplike space, we require that the space has two distinct basepoints.

Let $\Gamma$ be a space and $\gamma_1,\gamma_2\in \Gamma$ satisfying the following conditions:

(I) there exists a map $\nu :(\Gamma,\gamma_1,\gamma_2)\to(\Gamma\times \{\gamma_1\}\cup \{\gamma_2\}\times \Gamma, (\gamma_1,\gamma_1), (\gamma_2,\gamma_2))$
such that $\mbox{proj}_i \circ \nu \simeq\mbox{id}$ as maps of triples for each $i=1,2$,
where $\mbox{proj}_1,\mbox{proj}_2:(\Gamma\times \{\gamma_1\}\cup \{\gamma_2\}\times \Gamma, (\gamma_1,\gamma_1), (\gamma_2,\gamma_2)) \to (\Gamma,\gamma_1,\gamma_2)$
are the canonical projections;

(II) there exists a map $\eta :\Gamma\to \Gamma$ such that:

(a) $\eta(\gamma_1)=\gamma_2, \eta(\gamma_2)=\gamma_1$;

(b) $\nabla \circ (\overline {\mbox{id}} \vee \overline {\eta}) \circ \nu$ is homotopic to the constant map at $\gamma_1$,
 where
$$\overline{\mbox{id}} \vee \overline {\eta}:\Gamma\times \{\gamma_1\}\cup \{\gamma_2\}\times \Gamma, (\gamma_1,\gamma_1), (\gamma_2,\gamma_2)
\to \Gamma\times \{\gamma_1\}\cup \{\gamma_2\}\times \Gamma, (\gamma_1,\gamma_2), (\gamma_2,\gamma_1)
$$ with
$\overline{\mbox{id}}(\gamma,\gamma_1)=(\gamma,\gamma_2)$, $\overline {\eta}(\gamma_2,\gamma)=
(\gamma_2,\eta(\gamma))$ for $\gamma\in\Gamma$ and $\nabla :(\Gamma\times \{\gamma_2\}\cup \{\gamma_2\}\times \Gamma,
(\gamma_1,\gamma_2), (\gamma_2,\gamma_1)) \to (\Gamma,\gamma_1,\gamma_2)$
is the folding map;

(c) similarly, $\nabla \circ (\tilde {\mbox{id}} \vee \tilde {\eta}) \circ \nu$ is homotopic to the constant map at $\gamma_2$, where
$$\tilde{\mbox{id}} \vee \tilde {\eta}:\Gamma\times \{\gamma_1\}\cup \{\gamma_2\}\times \Gamma, (\gamma_1,\gamma_1),
(\gamma_2,\gamma_2) \to \Gamma\times \{\gamma_1\}\cup \{\gamma_2\}
\times \Gamma, (\gamma_2,\gamma_1), (\gamma_1,\gamma_2)
$$ with
$\tilde{\mbox{id}}(\gamma_2,\gamma)=(\gamma_1,\gamma)$,
$\tilde {\eta}(\gamma,\gamma_1)=((\eta(\gamma),(\gamma_1))$ for $\gamma\in\Gamma$;

(III) Moreover, we have co-associativity so that the following diagram
{\small
\begin{equation*}
\begin{CD}
(\Gamma,\gamma_1,\gamma_2)   @>{\nu}>>  (\Gamma\times \{\gamma_1\}\cup \{v_2\}\times \Gamma, (\gamma_1,\gamma_1), (v_2,v_2)) \\
@V{\nu}VV                            @VV{\tilde {\mbox{id}} \vee \tilde {\nu}}V  \\
(\Gamma\times \{\gamma_1\}\cup \{v_2\}\times \Gamma, (\gamma_1,\gamma_1), (v_2,v_2)) @>{\overline {\mbox{id}} \vee \overline {\nu}}>>   (\Gamma\times \{(\gamma_1,\gamma_1)\}\cup \{v_2\}\times (\Gamma\times \{\gamma_1\}\cup \{v_2\}\times \Gamma), \gamma_1^*, v_2^*)
\end{CD}
\end{equation*}
}
is commutative up to homotopy,
where $\gamma_1^*=(\gamma_1,(\gamma_1,\gamma_1)), \gamma_2^*=(\gamma_2,(\gamma_2,\gamma_2))$, and the maps
$\overline{\mbox{id}}(\gamma,\gamma_1)=(\gamma,(\gamma_1,\gamma_1))$, $\overline{\nu}(\gamma_2,\gamma)=
(\gamma_2,\nu(\gamma))$, $\tilde{\mbox{id}}(\gamma_2,\gamma)=((\gamma_2,\gamma_2),\gamma))$,
$\tilde{\nu}(\gamma,\gamma_1)=(\nu(\gamma),\gamma_1)$ for $\gamma\in\Gamma$.

Now, we generalize the notion of a co-grouplike space presented e.g., in \cite{oda-shimizu}.
A {\it co-grouplike space with two basepoints} $\Gamma = (\Gamma,\gamma_1,\gamma_2;\nu, \eta)$
consists of a topological space $\Gamma$ together with basepoints $\gamma_1,\gamma_2$ and maps
$\nu, \eta$ satisfying conditions (I) - (III).
For any space $W$, the {\it smash product} is given by
$$
\Gamma (W):=W\times \Gamma /\{(w,\gamma_1)\sim (w',\gamma_1), (w,\gamma_2)\sim (w',\gamma_2)\}
$$
for any $w,w'\in W$.

\par For instance, if $\Gamma =([0,1],0,1;\nu,\eta)$ with
$\nu(t)=\left\{\begin{array}{ll}
(2t,0)&\mbox{if}\;0\le t\le\frac{1}{2},\\
(1,2t-1)&\mbox{if}\;\frac{1}{2}\le t\le 1
\end{array}
\right.
$ and
$\eta(t)=1-t$ for $t\in[0,1]$ then $\Gamma (W)=\wh \Sigma W$, the un-reduced suspension of $W$.
\begin{remark}
Note that if $\gamma_1=\gamma_2$, we obtain the usual co-grouplike structure
and  $\Gamma_0:=\Gamma /\sim$ given by identifying the
basepoints $\gamma_1$ and $\gamma_2$ is a co-grouplike space as well.
\end{remark}

Next, we define the $W$-$\Gamma$-Rhodes groups.

Let $\Gamma$ be a co-grouplike space with two basepoints, $(X,G)$
a $G$-space and $W$ a space. The {\em $W$-$\Gamma$-Rhodes group} of $X$ with respect to $W$ is defined to be
$$
\sigma^{\Gamma}_W(X,x_0,G)=\{[f;g]|\;f:(\Gamma(W),\bar{\gamma}_1,\bar{\gamma}_2) \to (X,x_0,gx_0)\}.
$$
Write $\tau_W^{\Gamma_0}(X,x_0)$ for the $\Gamma_0$-$W$-Fox group
considered in \cite{ggw3}.

We can easily show:
\begin{proposition}{\em
Let $\pi: \sigma^{\Gamma}_W(X,x_0,G) \to G$ be the projection sending
$[f;g]\mapsto g$. By identifying the two basepoints of $\Gamma(W)$, the
quotient space $\Gamma (W)/\sim$ is canonically homeomorphic to
$\Gamma_0 \wedge (W\cup \{*\})$. Furthermore,
$$
\mbox{\em Ker}\, \pi \cong [\Gamma (W)/\sim, X] = \tau^{\Gamma_0}_W(X,x_0).
$$}
\end{proposition}

Then, we obtain a general $\Gamma$-Rhodes exact sequence, generalizing \eqref{generalized-rhodes-exact}.

\begin{theorem}\label{Gamma-fox-rhodes-split}
The following sequence
$$
1\to \tau^{\Gamma_0}_W(X,x_0) \to \sigma^{\Gamma}_W(X,x_0,G) \stackrel{\pi}{\to} G \to 1
$$
is exact.
\end{theorem}

We now derive the following generalized split exact sequence for the $W$-$\Gamma$-Rhodes groups.

\begin{corollary}
{\em Let $W$ be a space with a basepoint $w_0$ and $\Gamma$ be a co-grouplike space with two basepoints. The following sequence
\begin{equation}\label{general-Gamma-rhodes-split}
1\to [\Gamma_0\wedge ((V\times W)/V),\Omega X] \to \sigma^{\Gamma}_{V\times W}(X,x_0,G) \stackrel{\dashleftarrow}{\to} \sigma^{\Gamma}_V(X,x_0,G) \to 1
\end{equation}}
is split exact.
\end{corollary}
\begin{proof}
From Theorem \ref{Gamma-fox-rhodes-split}, we have the following two short exact sequences:
$$
1\to \tau^{\Gamma_0}_{V\times W}(X,x_0) \to \sigma^{\Gamma}_{V\times W}(X,x_0,G) \stackrel{\pi}{\to} G \to 1
$$
and
$$
1\to \tau^{\Gamma_0}_{V}(X,x_0) \to \sigma^{\Gamma}_{V}(X,x_0,G) \stackrel{\pi}{\to} G \to 1.
$$
Moreover, the following split exact sequence was shown in \cite[Theorem 4.1]{ggw3}:
$$
1\to [\Gamma_0\wedge ((V\times W)/V),\Omega X] \to \tau^{\Gamma_0}_{V\times W}(X,x_0,G) \stackrel{\dashleftarrow}{\to} \tau^{\Gamma_0}_V(X,x_0,G) \to 1.
$$
A straightforward diagram chasing argument involving these short exact sequences yields the desired split exact sequence.
\end{proof}

\section{Evaluation subgroups of $W$-$\Gamma$-Rhodes groups}

We end this note by extending a result concerning the evaluation subgroups of the Rhodes groups and the Fox torus homotopy groups obtained in \cite[Theorem 2.2]{ggw2}.

Given a $G$-space $X$, the function space $X^X$ is also a $G$-space where the action is pointwise, that is, $(gf)(x)=gf(x)$ for $f\in X^X$, $g\in G$ and $x\in X$. Let $\Gamma$ be a co-grouplike space with two basepoints and $W$ be a space.
\par The {\em evaluation subgroup} of the
$W$-$\Gamma$-Rhodes group of $X$ is defined by
$$
\mathcal G \sigma^{\Gamma}_W(X,x_0,G):={\rm Im}(ev_*:\sigma^{\Gamma}_W(X^X,\mbox{id}_X,G) \to \sigma^{\Gamma}_W(X,x_0,G)).
$$
Similarly, the {\em evaluation subgroup} of $\tau^{\Gamma_0}_W(X,x_0)$ is defined by
$$
\mathcal G \tau^{\Gamma_0}_W(X,x_0):={\rm Im}(ev_*:\tau^{\Gamma_0}_W(X^X,\mbox{id}_X) \to \tau^{\Gamma_0}_W(X,x_0)).
$$

It is straightforward to see that the proof of \cite[Theorem 2.2]{ggw2} is also valid in the setting of
$W$-$\Gamma$-Rhodes groups. Therefore, we have the following generalization.

\begin{theorem}
Let $G_0$ be the subgroup of $G$ consisting of elements $g$ considered as homeo\-morphisms of $X$ which are freely homotopic to the identity map $\mbox{\em id}_X$. Then the following sequence
$$
1\to \mathcal G \tau^{\Gamma_0}_W(X,x_0) \to \mathcal G \sigma^{\Gamma}_W(X,x_0,G) \to G_0 \to 1.
$$
is exact.
\end{theorem}

\end{document}